%%%%%%%%%%%%%%%%%%%%%% PLAIN_TEX %%%%%%%%%%%%%%%%%%%%%%%
\magnification 1200
\font\abs=cmr8
\font\refit=cmti8
\font\refbf=cmbx8

\font\ccc=cmcsc10

\def\tens{\otimes}
\def\fraz#1#2{{\strut\displaystyle #1\over\displaystyle #2}}

\def\esp#1{e^{\displaystyle#1}}
\def\ii#1{\item{\phantom{1}#1 .\phantom{x}}}
\def\jj#1{\item{#1 .\phantom{x}}}

\def\fze{{\cal F}_{\ell}(E(2))}

\def\uze{{\cal U}_{\ell}(E(2))}

\def\pro #1#2 {{\ccc (#1.#2) Proposition.}\phantom{X}}
\def\dfn #1#2 {{\ccc (#1.#2) Definition.}\phantom{X}}
\def\cor #1#2 {{\ccc (#1.#2) Corollary.}\phantom{X}}
\def\lem #1#2 {{\ccc (#1.#2) Lemma.}\phantom{X}}
\def\rem #1#2 {{\ccc (#1.#2) Remarks.}\phantom{X}}
\def\rmk #1#2 {{\ccc (#1.#2) Remark.}\phantom{X}}
\def\thm #1#2 {{\ccc (#1.#2) Theorem.}\phantom{X}}
\def\dim {{\sl Proof.}\phantom{X}}
\def\fidi{\hskip5pt \vrule height4pt width4pt depth0pt}

\def\et#1{e^{{#1}\,\theta}}
\def\ep#1{e^{{#1}\,P_2/2}}
\def\epp#1{e^{{#1}\,P_2}}
\def\xb{\bar x}
\def\chib{\bar \chi}
\def\po#1#2{(#1)_{\hskip-.5pt{\scriptscriptstyle #2}}}
\def\lj{\lambda({\cal J})}

\hsize= 15 truecm
\vsize= 22 truecm
\hoffset= 1. truecm
\voffset= 0.3 truecm

\baselineskip= 12 pt
\footline={\hss\tenrm\folio\hss} \pageno=1

%%%%%%%%%%%%%%%%%%% DOCUMENT %%%%%%%%%%%%%%%%%%%%%%%%
\vglue 3truecm
\centerline{\bf A NEW CLASS OF DEFORMED SPECIAL FUNCTIONS}
\smallskip
\centerline{\bf FROM QUANTUM HOMOGENEOUS SPACES.}
\bigskip
\bigskip
\centerline{{\it
F.Bonechi ${}^1$, R.Giachetti ${}^1$, M.A.del Olmo ${}^2$, 
E.Sorace ${}^1$ and M.Tarlini ${}^1$.}}
\bigskip
\baselineskip= 10 pt

{\hskip 0.7 truecm}${}^1${\abs Dipartimento
di Fisica, Universit\`a di Firenze  e INFN Sezione di Firenze, Italy.}

{\hskip 0.7 truecm}${}^2${\abs Departamento de F\'\i sica Te\'orica,
Universitad de Valladolid E-47011, Spain.}

\bigskip
\bigskip
{\refbf Abstract.} {\abs We study the most elementary aspects of
harmonic analysis on a homogeneous space of a deformation
of the two-dimensional Euclidean group, admitting generalizations
to dimensions three and four, whose quantum parameter has the
physical dimensions of a length. The homogeneous space is recognized
as a new quantum plane and the action of the Euclidean quantum group is
used to  determine an eigenvalue problem for the Casimir operator, that
constitutes the analogue of the Schr\"odinger equation in the presence of
such deformation. The solutions are given in the plane wave and in the
angular momentum bases and are expressed in terms of hypergeometric
series with non commuting parameters.
\hfil PACS: 02.20+b, 03.65.Fd, 11.30-j}

\bigskip
\bigskip
\bigskip

%%%%%%%%%%%%%%%%%%%% INTRODUCTION %%%%%%%%%%%%%%%%%%%%%%%
\baselineskip=14pt
\noindent {\bf 1. Introduction.}
\bigskip
Homogeneous spaces provide a unified framework for a wide class of
mathematical problems and often give a sound geometrical interpretation
to many results of classical and quantum physics.
The definition of special functions, the integral
transformations and the harmonic analysis are significant instances in
the context of pure mathematics; the classification of elementary
Hamiltonian systems by means of coadjoint
orbits and their quantization according to the Kirillov theory [1] is
one of the most important applications to physical problems.
The homogeneous spaces of kinematical groups, such as the Euclidean or
the Poincar\'e group, moreover, originate in a natural way
the fundamental wave equations of mathematical physics: indeed
these equations are determined by an invariant element of the
corresponding Lie algebra operating with a canonical action
on the functions on the homogeneous space. It appears therefore that
homogeneous spaces are a major constituent of the theory of Lie
groups and their applications.

As soon as the theory of quantum groups was founded, it seemed natural
to introduce the definition of quantum homogeneous spaces. Since,
in this case, one could not apply to the geometrical notion of
underlying manifold, the approach was necessarily algebraiuantum spheres were initially defined [2]; later on
a systematic work of generalization of the procedure and a detailed
study of possible applications to special functions was undertaken
[3-9]. One could observe that the largest efforts were devoted
to homogeneous spaces of compact quantum groups: no surprise in that,
since that part of the theory was best known.

As mentioned above, in applications of physical nature
a central position is occuped
by kinematical symmetries described, for instance, by the Heisenberg
and the Euclidean or pseudo-Euclidean groups.  A large number of papers
has dealt with the problem of defining the $q$-deformations of these
groups, starting with the Heisenberg [10-12] and with the Euclidean group
in two dimensions, $E(2)$ [10]. The situation for the latter has been
clarified in [13], where it has been
shown the possibility of two different deformations of $E(2)$.
The first one, that we denote by $E_q(2)$, has been the most studied
[14-23]: its enveloping algebra has the same relations
as in the classical case, while the coproducts of the `translation'
generators are twisted primitive with respect to the exponential of
the `rotation' generator. The second deformation, that we denote
by $E_\ell(2)$, obtained by a simple contraction procedure [24-26],
was later on generalized to dimensions three [27-28] and
four [29-30], producing, {\it e.g} the so called $\kappa$-Poincar\'e. Its
main feature is that the deformation parameter undergoes a
contraction  and acquires physical dimensions: by means of
an appropriate rescaling of the generators the parameter
can be reabsorbed, as was to be expected on a physical ground. We thus
get a deformation of $E(2)$ without parameter: the singular nature
of this deformation has been fully clarified in [13].

In a recent paper [31] we have studied the homogeneous spaces of 
$E_q(2)$. We have found two structures that were already introduced in 
literature  from different points of view and independently of the action
of $E_q(2)$, namely `quantum planes' [32,33] and `quantum
hyperboloids' [34]. Moreover, by quantizing some Poisson homogeneous
spaces, also a `quantum cylinder' has been obtained [35]. The existence
of the Haar functional has been proved for $E_q(2)$, [36,37]: by
projecting it on homogeneous spaces, the techniques of the usual harmonic
analysis can be extended to the realm of quantum homogeneous
spaces and the connection with $q$-special functions can be made
explicit. From a physical point of view, this makes conceivable the
study of the solutions of explicit models along the
usual lines of wave mechanics.

In this paper we define a new quantum plane as a quantum homogeneous
space of $E_\ell(2)$. We then specify the canonical action of
$E_\ell(2)$ on this space: according to the lines developed in [31], the
action will be used to define an eigenvalue equation for the Casimir of
$E_\ell(2)$, that constitutes a new deformed version of the free
Schr\"odinger equation. Due to the absence of the deformation parameter,
the results we are going to present are of a deeply different nature
with respect to those obtained by the similar analysis
developed in [31]: there, the Hahn-Exton functions and the
$q$-exponentials are recovered as $q$-deformation of the Bessel and
exponential functions. In the present case, the diagonalization of the
Casimir on the linear and angular momentum bases yields new special
functions that can be expressed in terms of hypergeometric series with
non commuting parameters. The reasons because we find it interesting to
present our results are twofold. In the first place this type of
harmonic analysis can turn out to be relevant, or even fundamental, for
the solution of any possible model presenting such a quantum group
symmetry. For instance, a different real form of this group
-- reproducing a 1+1 deformation of the Poincar\'e group --
has been fruitfully applied to phonon physics and has proved to be
the dynamical symmetry group for such a physical system on the lattice
[38]. Physical systems with similar properties have been studied
in [39,40]. Secondly,
since the deformation of  $E_\ell(2)$ has been extended to higher
dimensions, our results certainly provide a very useful support for
understanding the nature of the special functions on those extensions.

\bigskip
\bigskip

%%%%%%%%%%%%%%%%%%%%%% SECTION II %%%%%%%%%%%%%%%%%%%%%%
\noindent {\bf 2. Quantum homogeneous spaces of $E_\ell(2)\,$.}
\bigskip

In this section we define the homogeneous spaces of $E_\ell(2)$ that
will be used in the following analysis. In order to make the presentation
sufficiently self-consistent, we find it useful to recall in a sketchy
way the algebraic properties of $E_\ell(2)$, as obtained in [10] and
its duality with the quantized functions, found in [13].

\medskip
\dfn 21 The Hopf algebra generated by $\et{-i} ,\,a_1,\,a_2$, with
relations
$$[\et{-i},a_1]=\fraz{z}2\;(1-\et{-i})^2\,,\quad
[\et{-i},a_2]=i\,\fraz{z}2\;(\et{-2i}-1)\,,$$
$$[a_1,a_2]=i\,z\,a_1\ ,
$$
coalgebra operations
$$\eqalign{
   {}&\Delta(\et{-i}) = \et{-i}\tens \et{-i}\,,\cr
   {}&\Delta(a_1)\, = \,\cos(\theta)\tens
a_1-\sin(\theta)\tens a_2 +a_1\tens 1\,,\cr
   {}&\Delta(a_2)\, = \,\sin(\theta)\tens a_1 + \cos(\theta)\tens a_2
+a_2\tens 1\,, \cr}
$$
and antipode
$$\eqalign{
S(a_1)=& -\cos(\theta)\,a_1-\sin(\theta)\,a_2\,,\quad\quad
S(\et{-i})=\et{i}\,,\cr
S(a_2)=&\, \sin(\theta)\,a_1-\cos(\theta)\,a_2\,,\cr
}$$
will be called the {\it algebra of the quantized functions on $E(2)$}
and denoted by $\fze$.
\medskip
Assuming $z$ real, a compatible involution is given by
$$a_1^*=a_1\,,\quad a_2^*=a_2\,,\quad \theta^*=\theta\,.$$

The quantized enveloping algebra $\uze$ is generated by the unity and
the three elements $P_1\,,P_2\,,J$ satisfying
$$[J,P_1]=(i/z) \sinh(zP_2)\,,\quad [J,P_2]=-i\;P_1\,,\quad [P_1,P_2]=0\,,$$
and such that
\baselineskip=16pt
$$\Delta P_1 = \ep{-z}\tens P_1+P_1\tens \ep{z} ,\,\quad
\Delta P_2=P_2\tens 1+1\tens P_2\,,$$
$$\Delta J =\ep{-z}\tens J+J\tens \ep{z} ,\,$$
$$S(P_2)=-P_2\,,\quad\quad S(P_1)=-P_1\,,\quad\quad
S(J)=-J-i\,z\,P_1/2\ ,$$
\baselineskip=14pt
with vanishing counity and involution
$$J^*=J\,,\quad\quad P_1^*=P_1\,,\quad\quad P_2^*=P_2.$$

We finally write the duality pairing between $\uze$ and $\fze$,
\baselineskip=16pt
$$\eqalign{
\langle \nu_1,\theta^r a_1^s a_2^t\rangle
  = &
\delta_{r,0}\delta_{s,1}\delta_{t,0}\,,\quad
\langle \nu_2,\theta^r a_1^s a_2^t\rangle
  = \delta_{r,0}\delta_{s,0}\delta_{t,1}\,,\cr
{} &\langle \tau,\theta^r a_1^s a_2^t\rangle = \,\delta_{r,1}\delta_{s,0}
\delta_{t,0}\,,\cr}
$$
\baselineskip=14pt
where
$$\tau=-i\,\ep{-z}\,(J-i(z/4)\,P_1),$$
and
$$\nu_1=-i\,\ep{-z}\,P_1\,,\quad ~~~~~~~~~\nu_2=-i\,P_2\,.$$
Observing that $\Delta
\tau = \esp{-iz\,\nu_2}\tens \tau + \tau\tens 1$ and using the
condition $\langle u^*,a\rangle=\overline{\langle u,(S(a))^*\rangle}$,
with $u\in\uze$ and $a\in\fze$, we have $\nu_1^*=-\nu_1$,
$\nu_2^*=-\nu_2$, $\tau^*=-\tau - i\,z\,\nu_1$.
Moreover, if we consider the rescaled variables
$z\,P_1\,,\,z\,P_2\in\uze$ and
$a_1/z,\;a_2/z\in\fze$, we see that the deformation parameter
is reabsorbed: this means that all the structures corresponding to
different values of  $z$ are isomorphic among themselves.
It is nevertheless useful to maintain
explicitly the  deformation parameter
in order to perform more easily the classical limit ($z\rightarrow 0$).
\medskip
We recall the general definition of the two different
left actions of an element $Y$ of the quantized enveloping algebra
on a quantized function $f$, namely
$$\eqalign{
  \ell(Y)f&=(id\tens Y)\circ\Delta f=\sum\limits_{(f)}~f_{(1)}\,
  \langle Y,f_{(2)}\rangle\,,\cr
 \lambda(Y)f&=(S(Y)\tens id)\circ\Delta f=\sum\limits_{(f)}~
  \langle S(Y),f_{(1)}\rangle \,f_{(2)}\,.\cr}$$
For later use (and with obvious notation)
we also recall that these actions have the following properties:
$$\ell(YZ)f=\ell(Y)\ell(Z)\,f\,,~~~~~~~\lambda(YZ)f
=\lambda(Y)\lambda(Z)\,f ,$$
and
$$\ell(Y)fg=\sum\limits_{(Y)}~\ell(Y_{(1)})f~\ell(Y_{(2)})g\,,~~~~~~~
  \lambda(Y)fg=\sum\limits_{(Y)}~\lambda(Y_{(2)})f~\lambda(Y_{(1)})g\,.$$

Following the theory developed in [41], based on the ``infinitesimal
invariance'' of the functions on the quantum homogeneous spaces, we can
use the approach of [31] to look for  quantum homogeneous spaces of
$E_\ell(2)$.
\medskip
\lem 22 {\it Define
$~X=J-i(z/4)\, P_1\,.~$
The linear span of $X$ constitutes a $(*\circ S)$-invariant
two-sided coideal of $\uze$, twisted primitive with respect to
$\ep{-z}$.}
\smallskip
\dim By a straightforward calculation we have
$$*\circ S\,(X)=-X\,,\quad {\rm and}\quad \Delta X=\ep{-z}\tens X +
X\tens\ep{z}\,.\fidi$$
\medskip
\pro 23 {\it Let
$\,\,x=a_1-i\,a_2\,$, $\xb=a_1+i\,a_2\,$.
Then $x^*=\xb$ and
$$[x,\xb]=-z\,(x+\xb)\,,$$
Moreover $x$ and $\xb$ generate the invariant
subalgebra and left coideal
$$B_X=\{f\in \fze|\ell(X)\,f=0\}\,.$$
They thus define a quantum homogeneous space whose coaction reads:
$$\delta x=\et{-i}\tens x + x\tens 1\,, \quad~~~~~ \delta\xb=\et{i}\tens\xb
+ \xb\tens 1\,.$$
}
\smallskip
\dim First we observe that $X=i\,\ep{z}\,\tau\,,$ so that the kernel of
$\ell(X)$ is the same as that of $\ell(\tau)$. Let us write
$f=\sum_{l,m,n} f_{lmn}\,\et{-il}a_1^m a_2^n\,$. Then
$$\eqalign{
\ell(\tau)\,f=&\sum_{l,m,n}f_{lmn}\,\Big(\ell(e^{-i\,z\,\nu_2})\,\et{-il}
\ell(\tau)a_1^m a_2^n + (\ell(\tau)\,\et{-il})a_1^m a_2^n\Big)\cr
=& -i\,\sum_{l,m,n}\,l\;f_{lmn}\;\et{-il}a_1^m a_2^n\,\cr }$$
that vanishes for $l=0$. Then the space $B_X$ is generated by $a_1$ and
$a_2$ or by $x$ and $\xb$. It is immediate to calculate the
relationships between $x$ and $\xb$ as well as the coactions.\fidi
\medskip
\pro 24 {\it  The action $\lambda$ on $\fze$ has the following form:
$$\eqalign{
{}&\lambda(P_1)\;\et{-il}a_1^m a_2^n =
-i\,m\,\et{-il}a_1^{m-1}(a_2+i\,z/2)^n\,,\cr
{}&\lambda(P_2)\;\et{-il}a_1^m a_2^n =-i\,n\,\et{-il}a_1^m
a_2^{n-1}\,,\cr 
{}&\lambda(J)\,\et{-il}a_1^m a_2^n =i\,\et{-il}
\Big(\big(il\,a_1^m + m\, a_1^{m-1}(a_2-\fraz i2(m-1/2)z)\big)
(a_2+i\,z/2)^n\cr
{}&\phantom{\lambda(J)\,\et{-il}a_1^m a_2^n =i\,\et{-il}}
+\fraz i{2z} a_1^{m+1}\big((a_2+i\,z/2)^n - (a_2-i\,3z/2)^n\big)
\Big)\,.
\cr}$$
In particular}
$$\eqalign{
{}&\lambda(X)\;\et{-il}a_1^m a_2^n =i\,\et{-il}
\Big(\big(il\,a_1^m + m\,a_1^{m-1}(a_2-imz/2)\big)(a_2+i\,z/2)^n\cr
{}&\phantom{\lambda(X)\;\et{-il}a_1^m a_2^n =i\,\et{-il}}
+\fraz i{2z}a_1^{m+1}\big((a_2+i\,z/2)^n -
(a_2-i\,3z/2)^n\big)\Big)\,.
\cr}$$
\smallskip
\dim A straightforward calculation.\fidi
\bigskip
\bigskip

%%%%%%%%%%%%%%%%% SECTION III %%%%%%%%%%%%%%%%%%%%%%%%%%
\noindent {\bf 3. Free $\ell$-Schr\"odinger equation.}
\bigskip
The natural $\ell$-analog of the free Schr\"odinger equation is obtained
by considering the canonical action of the Casimir of $\uze$ on the
homogeneous spaces so far determined.

The Casimir of $\uze$ reads
$${\cal C}=4H^+H^-=(4/{z^2})\; \sinh ^2((z/2)\,P_2)+P_1^2\,,$$
where the elements
$$H^+=\fraz 1{2z}(\epp{z}-1)-\fraz 12 i\,\ep{z}\,P_1\,, \quad
H^-=\fraz 1{2z}(1-\epp{-z})+\fraz 12 i\,\ep{-z}\,P_1$$
are the deformations of the holomorphic and antiholomorphic
operators $P_2/2\mp i P_1/2$.
The coproducts of $H^{\pm}$ are
$$\Delta H^+=1\tens H^+ + H^+\tens\epp{z}\,,\quad\quad
\Delta H^-=\epp{-z}\tens H^- + H^-\tens 1\,.$$
Thus the $z$-deformed free Schr\"odinger equation reads:
$$4\,\lambda(H^+H^-)\;\psi=E\;\psi\,.\eqno(3.1)$$

In the remaining part of this section we shall diagonalize the operator
on the right hand side of (3.1) in the `plane wave' and `angular momentum'
bases, in analogy to the usual procedure carried on in quantum mechanics.
For later convenience we recall  the definition of
the {\it Pochammer symbol}
$$\po an =\prod_{k=0}^{n-1}(a+k)\,$$
and of the classical generalized hypergeometric series
$${}_rF_s\left[\matrix{a_1&\cdots&a_r\cr
                       b_1&\cdots&b_s\cr}\;;\;\zeta\right]\,=\,
  \sum^{\infty}_{\ell=0}\fraz{\phantom{l!\;}\po{a_1}\ell\cdots\po{a_r}\ell}
              {\ell!\;\po{b_1}\ell\cdots\po{b_s}\ell}\;\zeta^\ell\ .
$$
We also find it useful to introduce the variables
$$\chi=x/z\,,~~~~~~~~~~\bar{\chi}=-\bar{x}/z\,.$$

\bigskip
\bigskip
\noindent {\sl $(i)$ The plane wave  states.}
\medskip
In analogy with the standard procedure of diagonalization of the
Schr\"odinger operator on the plane wave basis, it is natural to
consider $H^+$ and $H^-$ as the appropriate candidates to be
diagonalized. We first determine the action of these operators.
\medskip
\lem 32 {\it The following relations hold:}
$$\eqalign{
{}&\lambda(H^+)\,\po{\chib}n=0\,,~~~~~~~
\lambda(H^+)\,\po{\chi}n=-\,\fraz nz\,\po{\chi}{n-1}\,,\cr
{}&\lambda(H^+)\,\po{1-\chi}n=-\,\fraz nz\,\po{1-\chi}n/(1-\chi)\,,
\cr}$$
{\it and}
$$\eqalign{
{}&\lambda(H^-)\,\po{\chi}n=0\,,~~~~~~~
\lambda(H^-)\,\po{\chib}n=\,\fraz nz\,\po{\chib}n/(\chib)\,,\cr
{}&\lambda(H^-)\,\po{1-\chib}n=\fraz nz\, \po{1-\chib}{n-1}\,.\cr
}$$
\smallskip
\dim By direct computation.\fidi
\medskip
\noindent These relations suggest to look for plane wave states of the 
form
$$\psi_{h^+h^-}=\sum_{m,n} h_{mn}\,\po{\chi}n\,\po{1-\chib}m\,,$$
where $h_{mn}$ are to be found from the equation
$$\lambda(H^+)\,\psi_{_{h^+h^-}}=h^+\,\psi_{_{h^+h^-}}\,.$$
\medskip
\pro 33 {\it The coefficients $h_{mn}$ are determined up to a
multiplicative constant and have the following expression
$$h_{mn}=\fraz 1{m!n!}(-zh^+)^m\,(zh^-)^n\,.$$}
\smallskip
\dim Indeed, using lemma (3.2),
$$\eqalign{
\lambda(H^+)\,\po{\chi}m\po{1-\chib}n&=-\fraz mz\po{\chi}{m-1}
\po{1-\chib}n\,,\cr
\lambda(H^-)\,\po{\chi}m\po{1-\chib}n&=\fraz nz\po{\chi}m
\po{1-\chib}{n-1}\,,\cr} $$
so that
$$\eqalign{
\lambda(H^+)\,\psi_{_{h^+h^-}}&=-\fraz 1z \sum_{m,n} (m+1)\,h_{m+1,n}
\po{\chi}m\po{1-\chib}{n}\,,\cr
\lambda(H^-)\,\psi_{_{h^+h^-}}&=\fraz 1z \sum_{m,n} (n+1)\,h_{m,n+1}
\po{\chi}m\po{1-\chib}{n}\,,\cr} $$
We thus find the following recurrence relations
$$-\fraz 1z (m+1)h_{m+1,n}=h^+\,h_{mn}\,,\qquad
\fraz 1z (n+1)h_{m,n+1}=h^-\,h_{mn}\,,$$
whose solution is straightforward and yields the result.\fidi
\medskip
The Casimir is obviously diagonal on these states:
$$\lambda({\cal C})\,\psi_{_{h^+h^-}}=4\,\lambda(H^+H^-)\,\psi_{_{h^+h^-}}
=4\,h^+h^-\,\psi_{_{h^+h^-}}\,.$$
Finally, the eigenfunctions $\psi_{_{h^+h^-}}$ can be expressed
in the form of a hypergeometric series. Indeed
$$\eqalign{
\psi_{_{h^+h^-}}&=\sum^{\infty}_{m=0}\fraz{(-zh^+)^m}{m!}\po{\chi}m\,
\sum^{\infty}_{n=0}\fraz{(zh^-)^n}{n!}\po{1-\chib}n\cr
{}&=(1+zh^+)^{\displaystyle{-\chi}}(1-zh^-)^{\displaystyle{\chib-1}}\,
\cr} $$
and, according to the general definition, this can be written as
$$\psi_{_{h^+h^-}}={}_1F_0\left[\matrix{\chi\cr -\cr}\;;\;-z\,h^+\right]\;
         {}_1F_0\left[\matrix{1-\chib\cr -\cr}\;;\;z\,h^-\right]\ .
$$

In the classical limit $z\rightarrow 0$ we recover the usual plane waves
$e^{i\,(h^+x-h^-\xb)}$, as expected.

\bigskip
\noindent {\sl $(ii)$ The angular momentum states.}
\medskip
We now consider the diagonalization of (3.1) on a basis which realizes
the deformed counterpart of the angular momentum states.
The duality structure and the requirement of a correct behaviour under
the involution indicate that such a result can be obtained
by diagonalizing the operator
$${\cal J}=\ep{-z}(J-i(z/4)\, P_1).$$
Observing that ${\cal J}^*={\cal J}$, let us therefore discuss the
equation (3.1) together with
$$\lj\;\psi=r\;\psi\,.$$

We now prove the following two lemmas.
\medskip
\lem 34 {\it The polynomials
$\po{\chi}n$ and $\po{\chib}n$ are eigenstates of $\lj$ with
eigenvalues $-n$ and $n$ respectively.}
\smallskip
\dim We prove the lemma by induction. For $n=1$ we have
$$\lj(\chi)=-\chi\,,\quad
\lj(\chib)=\chib\,,$$
as $\lj a_1=i\,a_2\,,$ and
$\lj a_2=-i\,a_1\,$. Assuming that
$$\lj\po{\chi}{n-1}=-(n-1)\po{\chi}{n-1}\,,\qquad
\lj\po{\chib}{n-1}=(n-1)\po{\chib}{n-1}\,,$$
we find
$$\lj\po{\chi}n=\lj\po{\chi}{n-1}(\chi +n-1)=
-n\,\po{\chi}n\,,$$
where we have used
$\Delta{\cal J}=\epp{-z}\tens {\cal J}+{\cal J}\tens 1$
and $\lambda(\epp{-z})\, \chi= \chi+1$. In similar way we get
$$\lj\po{\chib}n=n\,\po{\chib}n$$
observing that $\lambda(\epp{-z})\, \chib= \chib+1\,$.\fidi

\medskip
\lem 35 {\it The polynomial
$\rho_n=\po{\chib}n\,\po{1-\chi}n=\po{\chi}n\,\po{1-\chib}n$
is invariant under the action of $\lj$, {\it i.e.}
$\lj\,\rho_n=0\,.$
Moreover it can be written as $\rho_n=\rho(\rho+2)(\rho+6)\cdots
(\rho+n(n-1))\,$, where $\rho=\chib(1-\chi)\,$.}
\smallskip
\dim We start by proving the second part of the lemma. Also in this
case we adopt the induction technique. We have $\rho_1=\rho$ and
$\rho_n=\po{\chib}{n-1}\,(\chib +n-1)\,\po{1-\chi}{n-1}\,(n-\chi)\,$.
Since
the relation $(\chi-\chib)P(\chi)=P(\chi+1)(\chi-\chib)$ holds for any
polynomial
$P(\chi)$, it is straightforward to commute $(\chib+n-1)$ with
$\po{1-\chi}{n-1}$. As a result of the commutation we find
$$\rho_n=\po{\chib}{n-1}\,\po{1-\chi}{n-1}(\chib\,(1-\chi)+n(n-1))=
\rho_{n-1}(\rho + n(n-1))$$
and then $\rho_n=\rho(\rho+2)(\rho+6)\cdots(\rho+n(n-1))\,$.

The proof of the first statement is a direct consequence of the previous
result. Indeed from
$(\chib + \alpha)\,(1-\chi - \alpha)=(\chi+\alpha)\,(1-\chib-\alpha)\,,$
where $\alpha$ is any number, we find
$\po{\chib}n\,\po{1-\chi}n=\po{\chi}n\,\po{1-\chib}n\,.$
Moreover
$$\lj\,\rho_n=\lj\,\rho_{n-1}(\rho+n(n-1))
= \lj\rho_{n-1}\lambda(\epp{-z})(\rho+n(n-1))\,.$$
Due to the fact that  $\lj\,\rho=0$, we see that $\lj\,\rho_{n-1}=0$
implies $\lj\,\rho_n=0\,$.
The lemma is therefore proved.\fidi
\medskip
We shall write the eigenstates of $\lj$ as
$$\psi_r=\sum_{\ell}c^\ell_r\,\rho_\ell\,\po{\chib}r\,,\qquad
\psi_{-r}=\sum_{\ell}c^\ell_{-r}\,\rho_\ell\,\po{\chi}r\,.$$
 As a matter of fact, by a simple computation, we find
$$\lj\,\psi_r=r\,\psi_r\,,\qquad
\lj\,\psi_{-r}=-\,r\,\psi_{-r}\,.$$

The coefficients $c^\ell_r$ of the expansion will be determined
by using (3.1). It will appear that the 
choice of the polynomials $\rho_n$ proves to be essential when we try
to diagonalize the Casimir ${\cal C}=4\,H^+H^-$. To this purpose we find
it very useful to introduce an auxiliary element, which, according to the
following proposition, yields a very simple form for
the relations in $\uze$.
\medskip
\pro 36 {\it Let $U^\pm=(1/2)\,(1+\epp{\mp z}\mp iz\,\ep{\mp z}P_1)\,$ and
${\cal H}^\pm=U^\pm H^\pm\,$. Then
$$[{\cal J},{\cal H}^+]={\cal H}^+\,,\qquad
[{\cal J},{\cal H}^-]=-{\cal H}^-\,,\qquad [{\cal H}^+,{\cal H}^-]=0\,.$$
Moreover $~{\tilde{\cal C}}={\cal H}^+{\cal H}^-=H^+H^-(1+z^2H^+H^-)\,.$}
\smallskip
\dim By direct computation.\fidi
\medskip
We will now discuss the action of $\lambda({\cal H}^+)$ and
$\lambda({\cal H}^-)$ on $\psi_r$ and $\psi_{-r}$. On the one hand,
the proposition (3.6) and the involution property $({\cal H}^+)^*=
{\cal H}^-$ permit to write
$$\lambda({\cal H}^+)\,\psi_{\pm r}=\epsilon\,\psi_{\pm r+1}\,,\qquad
\lambda({\cal H}^-)\,\psi_{\pm r}={\bar \epsilon}\,\psi_{\pm r-1}\,.
\eqno(3.7)$$
On the other hand, the left hand side of (3.7) can direcly computed. 

Indeed:
\medskip
\lem 38 {\it The following relations hold
$$\eqalign{
\lambda({\cal H}^+)\psi_{-r}&=-\fraz 1z\,\sum_\ell \big[
c^\ell_{-r}\,(\ell+r)+c^{\ell+1}_{-r}(\ell+1)(\ell+r)(\ell+r+1)\big]\,
\rho_\ell\,\po{\chi}{r-1}\,\,,\cr
\lambda({\cal H}^-)\psi_{-r}&=-\fraz 1z\,\sum_\ell
c^{\ell+1}_{-r}\,(\ell+1)\,\rho_\ell\,\po{\chi}{r+1}.\cr}$$}
\smallskip
\dim
A straightforward calculation shows that
$$\lambda(U^+)\chi^n{\bar{\chi}}^m=\chi^n(\bar{\chi}+1)^m\,,\quad
\lambda(U^-){\bar{\chi}}^m \chi^n=(\bar{\chi}-1)^m \chi^n\,.$$
Using this result it is not difficult to complete the proof.
\fidi
\medskip
Equations (3.7) and (3.8) imply the recurrence relations
$$\eqalign{
\epsilon\,c^\ell_{-r+1} &= -\fraz 1z\big[c^\ell_{-r}\,(\ell+r)+
c^{\ell+1}_{-r}\,(\ell+1)(\ell+r)(\ell+r+1)\big]\,,\cr
{\bar
\epsilon}\,c^\ell_{-r-1} &= -\fraz 1z c^{\ell+1}_{-r}\,(\ell+1)\,.\cr}$$
It is easily verified that the coefficients
$$c^\ell_{-r}=(-kz/{\bar \epsilon})^r
\fraz{(kz^2)^\ell}{\ell!(\ell+r)!}\,,$$
solve the two recurrence relations, provided that $k$ is related to
$\epsilon{\bar \epsilon}$ by
$$|\epsilon|^2=\epsilon{\bar \epsilon}=k(1+z^2\,k)\,.$$
Recalling that ${\tilde {\cal C}}={\cal H}^+{\cal H}^-=
H^+H^-(1+z^2H^+H^-)\,$, we get
$$\lambda(H^+H^-)\,\psi_{-r}=k\,\psi_{-r}\,.$$
We have therefore proved the following
\medskip
\pro 39 {\it The states that diagonalize $\lj$ and $\lambda(H^+H^-)$
are
$$\psi_{-r}=\sum^{\infty}_{\ell=0}c^\ell_{-r}\,\rho_\ell\,\po{\chi}r\,,\qquad
\psi_{r}=\sum^{\infty}_{\ell=0}c^\ell_r\,\rho_\ell\,\po{\chib}r\,,$$
where
$$c^\ell_{-r}=(-kz/{\bar \epsilon})^r
\fraz{(kz^2)^\ell}{\ell!(\ell+r)!}\,,\qquad
c^\ell_r=(-kz/\epsilon)^r
\fraz{(kz^2)^\ell}{\ell!(\ell+r)!}\,.$$}
\medskip
Some final remarks are in order. We observe that
also in the angular momentum basis the states $\psi_{\pm r}$ can be
written, as {\it classical} hypergeometric series with
noncommutative coefficients. Indeed, by making
explicit the form of $\rho_\ell$ we get:
$$\eqalign{
\psi_{-r}=&\sum^{\infty}_{\ell=0}(-kz/{\bar \epsilon})^r
\fraz{(kz^2)^\ell}{\ell!(\ell+r)!}\po{\chib}\ell\,\po{1-\chi}\ell\,
         \po{\chi}r\cr
         =&\;{}_2F_1\left[\matrix{ \chib \qquad  1-\chi \cr
                           \phantom{\chib} r+1 \phantom{1-\chi}\cr}
             \;;\;kz^2\right]\;
           (-kz/{\bar \epsilon})^r\; \fraz{\po{\chi}r}{r!}\ ,
}$$
and
$$\psi_{r}=\;{}_2F_1\left[\matrix{ \chib \qquad  1-\chi \cr
                           \phantom{\chib} r+1 \phantom{1-\chi}\cr}
             \;;\;kz^2\right]\;
           (-kz/{\epsilon})^r\; \fraz{\po{\chib}r}{r!}\ .
$$
We want to stress that the situation is very different from the
deformation treated in [31]. The {\it quantum} deformation
is signified through the noncommutative variables $\chi$ and $\chib$
that appear in the coefficients $a_1$ and $a_2$ of ${}_2F_1$. The
coefficient $b_1$ and the hypergeometric variable $\zeta$ are, instead,
numbers. The classical limit $z\rightarrow 0$ of the hypergeometric
${}_2F_1$ yield again the usual Bessel functions in the commutative
variable $\bar x x$. This is not surprising. It is nevertheless rather
peculiar that this result is due to a kind of confluence phenomenon
caused by the non commutativity of the arguments.

\bigskip\bigskip

\noindent {\bf Acknowledgments.} M.A. O. expresses his gratitude for the 
hospitality during his stay in Firenze. This work has been partially 
supported by an Accion Integrada Espa\~na--Italia and by a DGICYT grant
(PB92--0255) from Spain.

%%%%%%%%%%%%%%%%%%% REFERENCES %%%%%%%%%%%%%%%%%%%%%%%
\baselineskip= 12 pt
\bigskip
\bigskip

\centerline{{\bf References.}}

\bigskip
\baselineskip= 10 pt
{\abs
\ii 1 Kirillov A.A., ``{\refit Elements of the Theory of Representations}'',
      (Springer Verlag, Berlin, 1990).
\smallskip
\ii 2 Podle\'s P., {\refit Lett. Math. Phys.}, {\refbf 14} (1987) 193.
\smallskip
\ii 3 Vaksman L.L. and Soibelman Y.S., {\refit Funct. Anal. Appl.},
     {\refbf 22} (1988) 170.
\smallskip
\ii 4 Masuda T., Mimachi K., Nakagami Y., Noumi M., Ueno K., {\refit
      J. Funct. Anal.}, {\refbf 99} (1991) 127.
\smallskip
\ii 5
      Noumi M., Mimachi K., ``{\refit Askey-Wilson polynomials as spherical
      functions on SU${}_q${\abs (2)}}'', in Lecture Notes in Mathematics
      n. 1510, 221, (Kulish P.P. ed., Springer--Verlag, 1992).
\smallskip
\ii 6 Koornwinder T.H., {\refit Proc. Kon. Ned. Akad. Wet.}, Series A,
      {\refbf 92} (1989) 97.
\smallskip
\ii 7
      Koelink H.T., Koornwinder T.H., {\refit Proc. Kon. Ned. Akad. Wet.},
      Series A, {\refbf 92} (1989) 443.
\smallskip
\ii 8
      Vilenkin N.Ja. and Klimyk A.U., ``{\refit Representation of Lie
      groups and Special Functions}'', Vol. 3, (Kluwer Acad. Publ.,
      Dordrecht, 1992).
\smallskip
\ii 9
      Noumi M., Yamada H., Mimachi K., ``Finite-dimensional representations
      of the quantum group {\refit GL}${}_q$({\refit n},{\refbf C}) and the
      zonal spherical functions on
     {\refit U}${}_q$({\refit n}--1)/{\refit U}${}_q$({\refit n})'',
     {\refit Japanese J. Math.} (to appear).
\smallskip
\jj {10}
      Celeghini E., Giachetti R., Sorace E. and Tarlini M., {\refit
      J. Math. Phys.}, {\refbf 31} (1990)  2548 .
\smallskip
\jj {11}
      Celeghini E., Giachetti R., Sorace E. and Tarlini M., {\refit
      J. Math. Phys.}, {\refbf 32} (1991) 1155.
\smallskip
\jj {12}
      Bonechi F., Giachetti R., Sorace E. and Tarlini M.,{\refit
      Commun. Math. Phys}, {\refbf 169} (1995) 463.
\smallskip
\jj {13}
      Ballesteros A., Celeghini E., Giachetti R., Sorace E. and Tarlini M.
      {\refit J. Phys. A: Math. Gen.}, {\refbf 26} (1993) 7495.
\smallskip
\jj {14}
      Vaksman L.L. and Korogodski L.I., {\refit Sov. Math. Dokl.},
      {\refbf 39} (1989) 173.
\smallskip
\jj {15}
      Koelink H.T., ``{\refit On quantum groups and q-special functions}'',
      Thesis (Leiden University, 1991) and {\refit Duke Math. J.},
      {\refbf 76} (1994) 483.
\smallskip
\jj {16}
      Baaj S., {\refit C.R. Acad. Sci. Paris}, {\refbf 314} (1992) 1021.
\smallskip
\jj {17}
      Woronowicz S.L., {\refit Lett. Math. Phys.}, {\refbf 23} (1991) 251.
\smallskip
\jj {18}
      Woronowicz S.L.,
      {\refit Commun. Math. Phys.}, {\refbf 144} (1992) 417.
\smallskip
\jj {19}
      Woronowicz S.L.,
      {\refit Commun. Math. Phys.}, {\refbf 149} (1992) 637.
\smallskip
\jj {20}
      Schupp P., Watts P. and Zumino B., {\refit Lett. Math. Phys.},
      {\refbf 24} (1992)  141.
\smallskip
\jj {21}
      Bonechi F., Celeghini E., Giachetti R., Pere\~na C.M., Sorace E. and
      Tarlini M., {\refit J. Phys. A: Math. Gen.}, {\refbf 27} (1994) 1307.
\smallskip
\jj {22}
      Ciccoli N. and Giachetti R.,
      {\refit Lett. Math. Phys.} {\refbf 34} (1995) 37.
\smallskip
\jj {23}
      Dabrowski L. and Sobczyk J., ``Left regular representation and
      contraction of {\refit sl}${}_q$(2) to {\refit e}${}_q$(2)'',
      Preprint IFT-864/94.
\smallskip
\jj {24}
      Celeghini E., Giachetti R., Sorace E. and Tarlini M.,
      ``{\refit Contractions of quantum groups}'', in Lecture Notes in
      Mathematics
      n. 1510, 221, (Kulish P.P. ed., Springer--Verlag, 1992).
\smallskip
\jj {25}
      Ballesteros A., Herranz F.J., del Olmo M.A. and Santander M.,
      {\refit J. Phys. A: Math. Gen.}, {\refbf 26} (1993) 5801.
\smallskip
\jj {26}
      Ballesteros A., Herranz F.J., del Olmo M.A. and Santander M.,
      {\refit J. Math. Phys.}, {\refbf 36} (1995) 631.
\smallskip
\jj {27}
      Celeghini E., Giachetti R., Sorace E. and Tarlini M., {\refit
      J. Math. Phys.}, {\refbf 32} (1991) 1159.
\smallskip
\jj {28}
      Ballesteros A., Herranz F.J., del Olmo M.A. and Santander M.,
      {\refit J. Phys. A: Math. Gen.}, {\refbf 27} (1994) 1283.
\smallskip
\jj {29}
      Lukierski J., Ruegg H., Nowicki A. and Tolstoy V.N. {\refit
      Phys. Lett. B}, {\refbf 264} (1991) 331.
\smallskip
\jj {30}
      Ballesteros A., Herranz F.J., del Olmo M.A. and Santander M.,
      {\refit J. Math. Phys.}, {\refbf 35} (1994) 4928.
\smallskip
\jj {31}
      Bonechi F., Ciccoli N., Giachetti R., Sorace E. and Tarlini M.,
      {\refit Commun. Math. Phys.}, {\refbf 175} (1996) 161.
\smallskip
\jj {32}
      Manin Yu.I., ``{\refit Quantum groups and noncommutative geometry}''
      in Publications of C.R.M. 1561 (University of Montreal, 1988).
\smallskip
\jj {33}
      Rieffel M.A., ``{\refit Deformation Quantization for actions of
      {\refit R}${}^d$}, Memoirs A.M.S., 506 (1993).
\smallskip
\jj {34}
      Schm\"udgen K., {\refit Commun. Math. Phys.}, {\refbf 159} (1994) 159.
\smallskip
\jj {35}
      Ciccoli N. ``{\refit Quantum Planes and Quantum Cylinders from
      Poisson Homogeneous Spaces}'', {\refit J. Phys. A: Math. and Gen.},
      in press.
\smallskip
\jj {36}
      Koelink H.T., ``{\refit On quantum groups and q-special functions}'',
      Thesis (Leiden University, 1991) and {\refit Duke Math. J.},
      {\refbf 76} (1994) 483.
\smallskip
\jj {37}
      Baaj S., {\refit C.R. Acad. Sci. Paris}, {\refbf 314} (1992) 1021.
\smallskip
\jj {38}
      Bonechi F., Celeghini E., Giachetti R., Sorace E. and Tarlini M.,
      {\refit Phys. Rev. Lett.}, {\refbf 68} (1992) 3718.
\smallskip
\jj {39}
      Bonechi F., Celeghini E., Giachetti R., Sorace E. and Tarlini M.,
      {\refit Phys. Rev. B},  {\refbf 32} (1992) 5727.
\smallskip
\jj {40}
      Bonechi F., Celeghini E., Giachetti R., Sorace E. and Tarlini M.,
      {\refit J. Phys. A}, {\refbf 25} (1992) L939.
\smallskip
\jj {41}
      Dijkhuizen M.S., Koornwinder T.H., {\refit Geom. Dedicata.},
      {\refbf 52} (1994) 291.

}

\bye